\documentstyle[11pt]{article}
\input amssym.def
\input amssym.tex

\newtheorem{theo}{Theorem}

\newtheorem{exam}{Example}

\newcommand{\eeq}{\end{equation}}
\newcommand{\beql}[1]{\begin{equation}\label{#1}}
\newcommand{\eqn}[1]{(\ref{#1})}
\newcommand{\hsp}{\hspace*{\parindent}}

\newcommand{\pf}{\noindent{\bf Proof.~}}
\newcommand{\dd}{\ldots}
\newcommand{\la}{\lambda}
\newcommand{\af}{\alpha}
\newcommand{\sR}{{\cal R}}
\newcommand{\sG}{{\cal G}}
\newcommand{\sH}{{\cal H}}

\newcommand{\sA}{{\cal A}}
\newcommand{\sB}{{\cal B}}
\newcommand{\sJ}{{\cal J}}
\newcommand{\sT}{{\cal T}}
\newcommand{\ZZ}{{\Bbb Z}}

\font\myname=msbm10 at 12pt
\newcommand\Bdbl[1]{\mbox{\myname #1}}
\makeatletter
\def\@sect#1#2#3#4#5#6[#7]#8{\ifnum #2>\c@secnumdepth
     \def\@svsec{}\else 
     \refstepcounter{#1}\edef\@svsec{\csname the#1\endcsname.\hskip .75em }\fi
     \@tempskipa #5\relax
      \ifdim \@tempskipa>\z@ 
        \begingroup #6\relax
          \@hangfrom{\hskip #3\relax\@svsec}{\interlinepenalty \@M #8\par}%
        \endgroup
       \csname #1mark\endcsname{#7}\addcontentsline
         {toc}{#1}{\ifnum #2>\c@secnumdepth \else
                      \protect\numberline{\csname the#1\endcsname}\fi
                    #7}\else
        \def\@svsechd{#6\hskip #3\@svsec #8\csname #1mark\endcsname
                      {#7}\addcontentsline
                           {toc}{#1}{\ifnum #2>\c@secnumdepth \else
                             \protect\numberline{\csname the#1\endcsname}\fi
                       #7}}\fi
     \@xsect{#5}}
\def\@begintheorem#1#2{\it \trivlist \item[\hskip \labelsep{\bf #1\ #2.}]}
\makeatother

\begin{document}
\begin{center}
{\Large {\bf Modular and $p$-adic cyclic codes*}} \\
\vspace{1\baselineskip}
{\em A. R. Calderbank} and {\em N. J. A. Sloane} \\
\vspace{.25\baselineskip}
Mathematical Sciences Research Center \\
AT\&T Bell Laboratories \\
Murray Hill, NJ 07974 \\
\vspace{1\baselineskip}
{\bf ABSTRACT}
\vspace{.5\baselineskip}
\end{center}
\setlength{\baselineskip}{1.5\baselineskip}

This paper presents some basic theorems giving the structure of
cyclic codes of length $n$ over the ring of
integers modulo $p^a$ and over the $p$-adic
numbers, where $p$ is a prime not dividing $n$.
An especially interesting example is the
2-adic cyclic code of length 7 with generator polynomial
$X^3+ \la X^2 + (\la -1) X-1$,
where $\la$ satisfies $\la^2 - \la +2 =0$.
This is the 2-adic generalization of both the binary
Hamming code and the quaternary octacode (the latter being equivalent
to the Nordstrom-Robinson code).
Other examples include the 2-adic Golay code of length 24 and the 3-adic
Golay code of length 12.
\footnotetext[1]{A version of this paper appeared in  {\em Designs, Codes and Cryptography},
{\bf 6} (1995), pp. 21--35. The references have now been updated.}

\section{Introduction}
\hsp
This paper was prompted by the following questions.
It is known \cite{Fo94}, \cite{HKCSS} that the binary
polynomial $X^3 + X+1$ that generates the cyclic Hamming code of length 7
lifts to a polynomial $X^3 + 2X^2 +X+3$ over $\ZZ_4$ that generates the octacode, equivalent
to the binary nonlinear Nordstrom-Robinson code.
What codes are obtained if we continue to lift this
polynomial to $\ZZ_8$, $\ZZ_{16} , \dd$, and even to the 2-adic integers
$\ZZ_{2^\infty}$?
What is the general structure of cyclic codes over these rings?
(Sol\'{e} \cite{So88} had already suggested in 1988 that $p$-adic
cyclic codes should be investigated.)

The answer to the first question is given in Example~\ref{ex1} of Section~4,
where we describe the ``2-adic Hamming code'' of
length 7 in detail.
This is in a certain sense the first interesting 2-adic code.
In Examples~\ref{ex2} and \ref{ex4} we give 2-adic
versions of the Golay code and more generally of extended quadratic
residue codes of length $8m$, where $8m-1$ is prime,
and a 3-adic version of the Golay code of length 12.
Furthermore, this Hamming code and the two
Golay codes (and more generally a large class of
quadratic residue codes) are all MDS codes.
In particular the 2-adic Golay code has minimal Hamming distance 13,
even though every projection of it onto the integers modulo $2^a$
has minimal distance 8.
Section~4 also gives $p$-adic generalizations for other
classical families of codes, including BCH, Reed-Muller and quadratic residue codes.

The answer to the second question is given in Theorems~\ref{Th5} and \ref{Th6} of Section~3,
which are the main theoretical results of this paper.
It will be seen that modular and $p$-adic cyclic codes have a simple
and elegant structure.

Although cyclic codes over the integers modulo $q$ have been
discussed by a number of authors
(\cite{Bl72}, \cite{Bl75}, \cite{CB93}, \cite{CS293}, \cite{RSpp}--\cite{Wa82}),
these results seem to have been overlooked.

The results in Section~3, although not at all obvious, are easily
verified by the methods of commutative algebra or
representation theory \cite{CR62}, \cite{ZS58}, so we shall mostly
not give proofs.

As far as we know, this paper is the first to consider $p$-adic
codes.
(However, several authors (\cite{AM69}, \cite{Ca80}, \cite{Ne82}) have studied ``global'' or
complex-valued codes in connection
with the representation theory of
$PSL_2 (n)$ and other groups, and our $p$-adic codes are analogues of those complex codes.)
For general background on $p$-adic numbers,
see \cite{Ba64}, \cite{BS66}, \cite{Go93}, \cite{Ko77}.
\section{Codes mod $p^a$ and $p$-adic codes}
\hsp
We use the symbol $\ZZ_{p^a}$ to denote the ring $\ZZ / p^a \ZZ$ of integers modulo $p^a$, for any prime $p$ and positive integer $a$, and $\ZZ_{p^\infty}$ for the ring of $p$-adic integers.
This slightly unconventional notation has the advantage of allowing us to
use $\ZZ_q$ (where $q=p^a$, $1 \le a \le \infty$) to denote any one of these rings,
and allows us to state our results in a uniform way.

An element $u \in \ZZ_{p^a}$ may be written uniquely as a finite sum
$$u= u_0 + pu_1 + p^2 u_2 + \cdots + p^{a-1} u_{a-1} ~,$$
and any element of $\ZZ_{p^\infty}$ as an infinite sum
$$u=u_0 + pu_1 + p^2 u_2 + \cdots ~,$$
where $0 \le u_i \le p-1$.
The units in $\ZZ_{p^a}$ or $\ZZ_{p^\infty}$ are precisely the $u$ for which
$u_0 \neq 0$.
$\ZZ_{p^a}$ has characteristic $p^a$, and $\ZZ_{p^\infty}$ has characteristic 0.

The following definitions and remarks are straightforward generalizations of notions for $\ZZ_4$ codes given in \cite{CS293} and \cite{HKCSS}.

Let $\ZZ_q = \ZZ_{p^a}$, where $1 \le a \le \infty$.
The set $\ZZ_q^n$ of $n$-tuples from $\ZZ_q$ is of course
a $\ZZ_q$-module, and by a {\em linear code} over $\ZZ_q$ we mean any
$\ZZ_q$ sub-module of $\ZZ_q^n$.
We equip $\ZZ_q^n$ with the inner product
$v \cdot w = v_1 w_1 + \cdots + v_n w_n$ evaluated in $\ZZ_q$, and define dual and self-dual codes in the usual way.

A nonzero linear code $C$ over $\ZZ_{p^a}$, for $a$ finite, has a generator matrix which after a suitable permutation of the coordinates can be
written in the form
\beql{ME1}
G = \left[
\matrix{I & A_{01} & A_{02} & A_{03} & \cdots & A_{0,a-1} & A_{0a} \cr
0 & pI & pA_{12} & pA_{13} & \cdots & pA_{1,a-1} & pA_{1a} \cr
0 & 0 & p^2 I & p^2 A_{23} & \cdots & p^2 A_{2,a-1} & p^2A_{2a} \cr
\cdot & \cdot & \cdot & \cdot & \cdots & \cdot & \cdot \cr
0 & 0 & 0 & 0 & \cdots & p^{a-1} I & p^{a-1} A_{a-1,a} \cr}
\right]~,
\eeq
where the columns are grouped into blocks of sizes
$k_0$, $k_1 , \dd , k_{a-1} , k_a$, and the $k_i$ are nonnegative integers adding to $n$.
This means that $C$ consists of all codewords
$$[v_0 ~~v_1 ~~v_2 ~~\cdots ~~ v_{a-1} ] G~,$$
where each $v_i$ is a vector of length $k_i$ with components
from $\ZZ_{p^{a-i}}$, so that $C$ contains $p^k$ codewords, where
$$k= \sum_{i=0}^{a-1} (a-i) k_i ~.$$
We say that $C$ has {\em type}\footnote{This definition of type differs from the one given in \cite{CS293} , \cite{HKCSS}.
The present definition has the advantage that it applies also to $p$-adic codes.}
\beql{ME5}
1^{k_0} p^{k_1} (p^2)^{k_2} \cdots (p^{a-1} )^{k_{a-1}} ~.
\eeq
The zero code (containing only the zero codeword) has type $1^0$.
It is easy to see that the code $C$ with generator matrix \eqn{ME1} has a dual $C^\perp$ with generator matrix of the form
\beql{ME3}
\left[
\matrix{
B_{0a} & B_{0,a-1} & \cdots & B_{03} & B_{02} & B_{01} & I \cr
pB_{1a} & pB_{1,a-1} & \cdots & pB_{13} & pB_{12} & pI & 0 \cr
p^2 B_{2a} & p^2 B_{2,a-1} & \cdots & p^2 B_{23} & p^2 I & 0 & 0 \cr
\cdot & \cdot & \cdots & \cdot & \cdot & \cdot & \cdot \cr
p^{a-1} B_{a-1,a} & p^{a-1} I & \cdots & 0 & 0 & 0 & 0 \cr}
\right] ~,
\eeq
where the column blocks have the same sizes as in \eqn{ME1}.
The dual code therefore contains $p^{k_\perp}$ codewords, where
$$k_\perp = \sum_{i=1}^a ik_i ~,$$
and has type
\beql{ME9}
1^{k_a} p^{k_{a-1}} (p^2)^{k_{a-2}} \cdots
(p^{a-1})^{k_1} ~.
\eeq
Also $|C| |C^\perp| = p^{k+k_\perp} = p^{an}$, and $(C^\perp)^\perp =C$.

Similarly, a nonzero linear code $C$ over $\ZZ_{p^\infty}$ has a generator matrix which can be written in the form
\beql{ME6}
G= \left[
\matrix{p^{m_0}I & p^{m_0}A_{01} & p^{m_0} A_{02} & \cdots &
p^{m_0} A_{0,b-1} & p^{m_0} A_{0,b} \cr
0 & p^{m_1} I & p^{m_2} A_{12} & \cdots & \cdot & \cdot \cr
\cdot & \cdot & \cdot & \cdots & \cdot & \cdot \cr
0 & 0 & 0 & \cdots & p^{m_{b-1}}I & p^{m_{b-1}} A_{b-1,b} \cr}
\right]~,
\eeq
where $0 \le m_0 < m_1 < \cdots  m_{b-1}$, for some integer $b$,
the column blocks have sizes $k_0$, $k_1 , \dd , k_b$ and the $k_i$ are nonnegative integers adding to $n$.
This means that $C$
consists of all codewords
$$[v_0 ~~v_1 ~~v_2 ~~\cdots ~ v_b ] G~,$$
where each $v_i$ is a vector of length $k_i$ with components from
$\ZZ_{p^\infty}$.
We say that $C$ has type
\beql{ME7}
(p^{m_0} )^{k_0} (p^{m_1})^{k_1} \cdots
(p^{m_{b-1}})^{k_{b-1}} ~.
\eeq
Now the code contains infinitely many codewords (although
it is still finitely generated).

If $m_0 > 0$ in \eqn{ME6}, all the codewords are multiples of $p^{m_0}$,
and (since $\ZZ_{p^\infty}$ has characteristic 0) we may divide the whole code by
$p^{m_0}$.
We shall therefore usually only consider codes in which
$m_0 =0$.
In this case the dual code has a generator matrix similar
to \eqn{ME3}, with type
\beql{ME8}
1^{k_b} (p^{m_1})^{k_{b-1}} \cdots
(p^{m_{b-1}} )^{k_1} ~,
\eeq
and $(C^\perp)^\perp =C$.
(If $m_0 > 0$ then $(C^\perp )^\perp = p^{- m_0} C$.)

The {\em automorphism group} $Aut (C)$ of a linear code $C$ over
$\ZZ_q$ is defined to be the set of all monomial matrices
over $\ZZ_q$ that preserve the code.
Since it contains all scalar matrices $uI$, where $u$ is a unit in $\ZZ_q$,
this group is infinite if $q=p^\infty$.
We therefore define the {\em projective automorphism group} to be
the quotient group $Aut(C)/ \{ uI :~u = \mbox{unit} \}$.

A {\em cyclic} code $C$ of length $n$ over $\ZZ_q$
$(q=p^a$, $1 \le a \le \infty)$ is a
linear code with the property that if
$(c_0 , c_1 , \dd , c_{n-1} ) \in C$ then
$(c_1 , c_2 , \dd , c_{n-1} , c_0) \in C$.
We assume throughout that $n$ and $p$ are relatively prime.
As usual we represent codewords by polynomials, so cyclic
codes are precisely the ideals in the ring
$$\sR = \ZZ_q [X] / (X^n -1) ~.$$
\section{Rings}
\hsp
We now discuss the properties of the ring $\sR$ and of certain
Galois rings $GR(q^m)$.

Let $q=p^a$ $(1 \le a \le \infty )$, and let $\pi_1 (X) \in \ZZ_p [X]$ be a monic primitive irreducible polynomial of degree $m$, so that $\pi_1 (X)$
divides $X^n -1$ mod$\,p$, where $n=p^m -1$.
The following are straightforward generalizations of results
given in \cite{HKCSS}, \cite{Mc74}, \cite{Ya90}.
There is a unique monic irreducible polynomial $\pi_a (X) \in \ZZ_q [X]$ such
that $\pi_a (X) \equiv \pi_1 (X)$ $\bmod\,p$ and $\pi_a (X)$ divides $X^n -1$ over $\ZZ_q$ (see Theorem~\ref{Th1} below).

Let $\xi$ be a root of $\pi_a (X)$, so that $\xi^n =1$.
Then the {\em Galois ring} $GR (q^m)$ is by definition the ring $\ZZ_q [ \xi ]$.
There are two canonical ways to represent the elements of this ring.
In the first representation, every element has a unique
expansion
$$u=u_0 + pu_1 + p^2 u_2 + \cdots + p^{a-1} u_{a-1}$$
(an infinite sum if $a= \infty$), where $u_i \in \sJ = \{0,1, \xi , \xi^2 , \dd , \xi^{n-1} \}$.
The map $\tau : u \mapsto u_0$ is given by
$$\tau (u) = u^{p^m} ~, ~~~
u \in \ZZ_q [ \xi ] ~,$$
and satisfies
$$\tau (uv) = \tau (u) \tau (v) ~, ~~~ u,v \in \ZZ_q [ \xi ] ~.$$
In the second representation $u$ is written as
$$u= \sum_{r=0}^{n-1} v_r \xi^r ~, ~~~v_r \in \ZZ_q ~.$$
The {\em Frobenius map} $\phi$ from $\ZZ_q [ \xi ]$ to $\ZZ_q [\xi ]$ takes
$$\sum_{r=0}^{a-1} p^r u_r ~~~\mbox{to}~~~
\sum_{r=0}^{a-1} p^r u_r^p ~.$$
Then $\phi$ generates the Galois group of $\ZZ_q [ \xi ]$ over
$\ZZ_q$, and $\phi^m$ is the identity map.

The following theorem plays a central role in studying cyclic codes over
$\ZZ_q$.
It shows that the irreducible factors of $X^n -1$ over $\ZZ_q$ are in one-to-one
correspondence with the factors over $\ZZ_p$.
\begin{theo}
\label{Th1}
Let $q=p^a$, $1 \le a \le \infty$.
If $h_1 (X) \in \ZZ_p [X]$ is a monic irreducible divisor of
$X^n -1$ over $\ZZ_p$, then there is a unique monic
irreducible polynomial $h_a (X) \in \ZZ_q [X]$ which
divides $X^n -1$ over $\ZZ_q$ and is congruent to $h_1 (X)$ mod$\,p$.
\end{theo}
\pf
This result can be obtained from Hensel's Lemma, but we prefer to sketch a constructive proof (by induction).

For $1 \le r < \infty$, suppose $h_r (X) \in \ZZ_{p^r} [X]$ is a monic
irreducible polynomial such that
$h_r (X) \equiv h_1 (X)$ $\bmod\, p$, and $h_r(X) ~|~ X^n -1$ over
$\ZZ_{p^r}$.
We will show that $h_r(X)$ can be lifted uniquely to a monic irreducible polynomial $h_{r+1} (X) \in \ZZ_{p^{r+1}} [X]$ which divides $X^n -1$
over $\ZZ_{p^{r+1}}$.
Then $h_\infty (X)$ is defined as the ($p$-adic) limit of
$h_r (X)$ as $r \to \infty$.

Let $h(X) \in \ZZ_{p^{r+1}} [X]$ be any lift of
$h_r (X)$, say
$h(X) = h_r (X) + p^r g (X)$, and let $\alpha$ be a root of
$h_r (X)$ and $\beta$ a corresponding root of $h(X)$, so that
$\beta = \alpha + p^r \delta$.
Then
\begin{eqnarray*}
\alpha^n & = & 1+ p^r \epsilon~,~~~
\beta^p ~=~ ( \alpha + p^r \delta )^p ~=~
\alpha^p ~, \\
\beta^{np} & = & (1+ p^r \epsilon)^p ~=~ 1 ~.
\end{eqnarray*}
Therefore the monic polynomial whose roots are the $p$-th powers of the roots of $h(X)$ divides $X^n -1$, and $\bmod\,p^r$ has the same roots as
$h_r (X)$, and so may be taken as $h_{r+1} (X)$.
This polynomial is irreducible since its roots form one orbit
under the Frobenius map.
To show that $h_{r+1} (X)$ is unique,
we argue as follows.
Let $h(X)$ and $h' (X)$ be two different possibilities for
$h_{r+1} (X)$, and let $\beta$ and $\gamma$ be zeros of $h$ and
$h'$ respectively, with $\beta \equiv \gamma$ $\bmod\,p^r$, say
$\beta = \gamma + p^r \delta$.
Then $\beta^n = \gamma^n =1$, $\beta^p = \gamma^p$, hence $(\beta / \gamma)^n = ( \beta / \gamma)^p =1$.
Since $n$ and $p$ are relatively prime,
$\beta = \gamma$, and so $h=h'$. \hfill $\blacksquare$

We now investigate the structure of ideals in $\sR$.
The units in $\sR$ are precisely the elements $u = \sum\limits_{r=0}^{n-1} u_r X^r$, $u_r \in \ZZ_q$, such that
at least one of the $u_r$ is a unit in $\ZZ_q$.
We denote the natural map from $\sR$ to
$\ZZ_p [X] / (X^n -1)$ by $\mu$.

If $\sA$ is an ideal in $\sR$ with generators $f_1 , f_2 , \dd$, we write
$\sA = (f_1 , f_2 , \dd )$.
The radical $Rad ( \sA )$ of $\sA$ is the set of all elements of $\sR$, some power of which is in $\sA$.
The radical of the ideal $\{0\}$ is called the radical of
$\sR$, and denoted by $Rad ( \sR )$.
Then $Rad ( \sR) = (p)$ if $q=p^a$ is finite, or (0) if $q=p^\infty$.

The ring $\ZZ_{p^\infty}$ is a principal ideal domain, hence Noetherian.
This implies that $\ZZ_{p^a} [X]$ and $\sR = \ZZ_{p^a} [X] / (X^n -1)$ are Noetherian for all $1 \le a \le \infty$.
$\sR$ satisfies the descending chain condition if $q=p^a$ is finite (since then $\sR$ is finite), but not if $q=p^\infty$ (we will see examples later).
Hence every maximal ideal in $\sR$ is prime, and if $q$ is finite every prime ideal different from (0) and (1) is maximal
(\cite{ZS58}, pp.~150, 203).

It is well-known that the prime ideals in
$\ZZ_p [X] / (X^n -1)$ are (0), (1) and $(\pi_1)$, where
$\pi_1$ is any monic irreducible divisor of $X^n -1$ over $\ZZ_p$.
\begin{theo}
\label{Th2}
If $q=p^a$ is finite the prime ideals in $\sR$ are $(0)$, $(1)$ and $(\pi_a , p)$, where $\pi_a$ is any monic irreducible divisor of $X^n -1$ over $\ZZ_q$.
If $q=p^\infty$ there are in addition the prime (but nonmaximal) ideals $(\pi_a)$.
\end{theo}
\pf
Let $\sA$ be a prime ideal in $\sR$ different from (0) and (1).
Then $\mu ( \sA) = ( \pi_1)$, say, so
$\sA$ contains $\pi_a$, where $\mu (\pi_a) = \pi_1$.
If $q$ is finite then $p \in \sA$, or else $\sR / \sA$ would contain zero divisors,
so $\sA \supset ( \pi_a , p)$, and it is easily seen that this ideal
is maximal.
If $q$ is infinite and $p \not\in \sA$ then the only other
possibility is $\sA = ( \pi_a)$. \hfill $\blacksquare$

Note that the ideal $(p)$ is not prime, since it contains
the product of all the $\pi_a$ --- which is 0 --- but
none of the $\pi_a$ themselves.

It is also known that every ideal $\sA$ in $\ZZ_p [X] / (X^n -1)$ contains an
idempotent $e_1$ (say), such that $\sA = (e_1)$ (\cite{MS77}, Chapter~8, Theorem~1; \cite{CR62}, \S24.2).
\begin{theo}
\label{Th3}
Every prime ideal $\sA = ( \pi_a , p)$ in $\sR$ contains an idempotent $e_a$ with $e_a^2 = e_a$,
$\sA = (e_a ,p)$.
Furthermore, if $q$ is infinite then every prime ideal $\sA = ( \pi_a)$ has an idempotent generator.
\end{theo}
\pf
We establish the first assertion by induction.
Let $(\pi_r ,p)$ be the projection of $\sA$ onto $\ZZ_{p^r} [X]/(X^n-1)$, and
suppose $e_r \in (\pi_r,p)$ is an idempotent with $(e_r,p) = (\pi_r, p)$.
Then
$e_r^2 = e_r + p^r h$ in $\ZZ_{p^{r+1}} [X] / (X^n-1)$,
for some $h$ in $\ZZ_{p^{r+1}} [X] / (X^n-1)$.
If we take $e_{r+1} = e_r + p^r \theta$, then $e_{r+1}^2 - e_{r+1} = p^r (h- \theta (1-2e_r))$,
and $e_{r+1}$ is an idempotent in $\ZZ_{p^{r+1}} [X]/ (X^n -1)$ if we choose
$\theta =h$ (if $p=2$) or
$\theta = h(1-2e_r)^{-1}$ (if $p > 2$).
(Note that $(1-2e_r)^2 = 1+ 4p^r h$, so $1-2e_r$ is a unit.)
It is easily verified that $(e_{r+1} , p) = (\pi_{r+1}, p)$.
By repeating this process we obtain an idempotent $e_a \in \sA$ with
$(e_a , p) = (\pi_a,p)$.

To prove the second assertion, since $\pi_a$ and $(X^n -1)/ \pi_a$ are relatively
prime, we can find $h \in \ZZ_{p^\infty} [X]$ such that
$$h \pi_a -1 \equiv 0 ~~~\bmod \, (X^n-1)/ \pi_a ~,$$
so $h \pi_a (h\pi_a -1) =0$ in $\sR$,
and $h \pi_a$ is the desired idempotent. \hfill $\blacksquare$

Next, every primary ideal is a power of a prime ideal.
\begin{theo}
\label{Th4}
The primary ideals in $\sR$ are $(0)$, $(1)$, $(\pi_a)$ and $(\pi_a , p^i)$, where $\pi_a$ is an irreducible divisor of $X^n -1$ over $\ZZ_q$ and $1 \le i < a$.
\end{theo}

We omit the proof.
The key steps are
(i)~to show that if $\sA = (\pi_a, p) = (e_a, p)$ is
a prime ideal then
\beql{ME10}
\sA^i = (\pi_a , p)^i = (\pi_a,p^i) = (e_a , p^i) ~,
\eeq
for $1 \le i < a$, and
(ii)~to show that if $\sB$ is a primary ideal whose
associated prime ideal is $\sA = (\pi_a, p)$ then
(by \cite{ZS58}, p.~200, Ex.~2) there is an integer $j$ such that
$\sA^j \subseteq \sB \subseteq \sA$, and from this
that $\sB = \sA^i$ for some $i$.

Note that when $q=p^a$ is finite then $(\pi_a ,p)^a = (\pi_a)$, and
$$(\pi_a ,p) \supset (\pi_a ,p^2) \supset \cdots \supset
(\pi_a ,p^{a-1}) \supset (\pi_a)$$
is a finite descending sequence.
When $q=p^\infty$, however,
$$(\pi_\infty ,p) \supset (\pi_\infty ,p^2) \supset (\pi_\infty,p^3) \supset \cdots \supset (\pi_\infty)$$
is an infinite descending sequence of primary ideals, the first and last
of which are prime.
In this case we adopt the convention that $(\pi_\infty ,p)^\infty$
denotes $(\pi_\infty )$.
\begin{theo}
\label{Th5}
Let $\pi_a^{(i)}$, $i=1, \dd , A$, denote the distinct
monic irreducible divisors of $X^n -1$ over $\ZZ_q$.
Any ideal in $\sR$ can be written in a unique way as
\beql{ME11}
\sA = \prod_{i=1}^A (\pi_a^{(i)} , p)^{m_i} ~,
\eeq
where $0 \le m_i \le a$.
In particular if $a$ is finite there are $(a+1)^A$ distinct ideals.
\end{theo}

This is a consequence of Theorem~\ref{Th4}
and the Lasker-Noether decomposition theorem (\cite{ZS58}, p.~209).
The product symbol in \eqn{ME11} may also be replaced by an
intersection symbol.
\begin{theo}
\label{Th6}
If $q=p^a$, $1 \le a < \infty$, any ideal in $\sR$ has the form
\beql{ME12}
(f_0 , pf_1, p^2 f_2 , \dd , p^{a-1} f_{a-1} ) ~,
\eeq
where the $f_i$ are divisors of $X^n -1$ satisfying
\beql{ME13}
f_{a-1} ~\Bigl|~ f_{a-2} ~\Bigl|~
\cdots ~\Bigl|~ f_1 ~\Bigl|~ f_0 ~.
\eeq
If $q=p^\infty$, any ideal in $\sR$ has the form
\beql{ME14}
(p^{m_0} f_0 , p^{m_1} f_1 , \dd , p^{m_{b-1}} f_{b-1} ) ~,
\eeq
where $0 \le m_0 < m_1 < \cdots < m_{b-1}$, for some
$b$, and
$$f_{b-1} ~\Bigl|~ f_{b-2} ~\Bigl|~ \cdots ~\Bigl|~ f_1 ~\Bigl|~ f_0 ~.$$
\end{theo}
\noindent
{\bf Proof.}
This follows by expanding the product in \eqn{ME11} and using
\eqn{ME10}. \hfill $\blacksquare$

\vspace*{+.1in}
\noindent{\bf Corollary.}
{\em Every ideal in $\sR$ is principal.}

\vspace*{+.1in}
\noindent
{\bf Proof.}
(i)~If $q=p^a$, $1 \le a < \infty$, then the ideal defined by
\eqn{ME12} has the generator
$$g=f_0 + pf_1 + p^2 f_2 + \cdots + p^{a-1} f_{a-1} ~.$$
We prove this for $a=2$ and 3, leaving the general case to the reader.
Let $\widehat{f}_0 = (X^n -1) / f_0$, $\widehat{f}_i = f_{i-1} / f_i$ for
$1 \le i < a$.
Case $a=2$:
Then $g = f_0 + pf_1$, and $(g)$ contains $pg = p f_0 = pf_1 \widehat{f}_1$ and
$\widehat{f}_0 g = pf_1 \widehat{f}_0$, hence
$pf_1$ (since $\widehat{f}_0$ and $\widehat{f}_1$ have no
common factors), hence $f_0$.
Case $a=3$:
Now $g= f_0 + pf_1 + p^2 f_2$, and $(g)$ contains $p^2 g = p^2 f_2 \widehat{f}_1 \widehat{f}_2$, $p \widehat{f}_0 g = p^2 f_2 \widehat{f}_0 \widehat{f}_2$, and
$\widehat{f}_0 \widehat{f}_1 g = p^2 f_2 \widehat{f}_0 \widehat{f}_1$,
hence $p^2 f_2$, hence
$f_0 + pf_1$.
So $(g) = (f_0 + pf_1 , p^2 f_2 )$.
Arguing as in case $a=2$ it follows that
$(g) = (f_0 , pf_1 , p^2 f_2)$.

(b)~Suppose $q= p^\infty$.
Let $g_a$ be a generator for the principal ideal given by the projection of the ideal onto $\ZZ_{2^a}$, for $a=1,2, \dd~$.
Since $\sR$ is compact in the $p$-adic metric, the sequence
$\{g_a\}$ has a subsequence which converges to a limit $g$ (say).
Then $g$ generates the ideal. \hfill $\blacksquare$

Finally, although we have not made any use of this,
it is worth noting that $\sR$ has a decomposition into a direct
product of Galois rings:
$$\ZZ_{p^a} [X] / (X^n -1) \cong \prod_{i=1}^A
\ZZ_{p^a} [X] / ( \pi_a^{(i)} )~.$$
\section{Generalizations of classical codes to $\Bdbl{Z}_q$}
\hsp
Theorem~\ref{Th1} provides a mechanism for generalizing any class of
cyclic codes from $GF(p)$ to $\ZZ_{p^a}$ (for finite $a$) and even to the
$p$-adic integers $\ZZ_{p^\infty}$.
For example we define a {\em BCH code} of
length $n$ over $\ZZ_q$ $(q=p^a$, $1 \le a \le \infty)$
to be the cyclic code whose generator polynomial is obtained by lifting
the generator polynomial for a BCH code over $GF(p)$ to
$\ZZ_q$.
The resulting polynomial has a string of consecutive roots in the appropriate Galois
ring $GF(q^m)$.
(For finite $q$ this is essentially the same as Shankar's \cite{Sh79}
definition of BCH codes over $\ZZ_q$.)
The code has type $1^k$, where $k$ is the dimension of the BCH
code over $GF(p)$.
One of the main unsolved questions here is to determine how the
minimal Lee distance of these BCH codes varies as $a \to \infty$.
(Similar questions can be asked about all the codes in this section.)
We investigate the first nontrivial case of these BCH codes later in this section.

We define {\em Reed-Muller codes} (since they are extended cyclic
codes \cite{AK92}, \cite{MS77}) and
{\em quadratic-residue codes} over $\ZZ_q$ in an
analogous way.

If $C$ is a code of length $n$ over $\ZZ_q$ with generator matrix
\eqn{ME1} or \eqn{ME6} and type \eqn{ME5} or \eqn{ME7}, we define $k$ by
$$
k = \sum_{i=0}^{a-1} k_i ~~\mbox{(for \protect\eqn{ME5})},~~
\sum_{i=0}^{b-1} k_i~~
\mbox{(for \protect\eqn{ME7})}~.
$$
The usual argument (\cite{MS77}, Chapter~2) then gives the
Singleton bound:
\beql{ME15}
d \le n-k+1 ~,
\eeq
where $d$ is the minimal Hamming distance of the code.
We say that $C$ is {\em maximal distance separable},
or MDS, if equality holds in \eqn{ME15}.
Since codes over $\ZZ_{p^\infty}$ have infinitely many codewords, it is
better to use the equivalent definition (see \cite{MS77}, Chapter~11, Corollary~3) that a code is
MDS if and only if every $k$ columns
of the generator matrix are linearly independent over $\ZZ_q$.
\begin{exam}
\label{ex1}
{\bf The 2-adic Hamming code of length 7}.
{\rm In the binary case,
$X^n -1$ factors trivially over $\ZZ_q$, $q=2^a$,
$1 \le a \le \infty$, for $n=1,3$ and 5.
The first nontrivial factorization is for $n=7$, where it is easy\footnote{Guided by the factorizations mod~2 and mod~4, one guesses that
$X^6 + X^5 + \cdots + 1 = (X^3 + \lambda X^2 + \mu X-1) \cdot$ reciprocal;
hence
$\mu = \lambda -1$, $\lambda^2 = \lambda -2$.} to find the
2-adic factorization
\beql{ME16}
X^7 - 1 = (X-1)(X^3 + \la X^2 + (\la-1)X-1)
(X^3 - (\la-1)X^2 - \la X -1)~,
\eeq
where
\beql{ME17}
\la = 0+2+4+32+128+256+ \cdots
\eeq
is a 2-adic number satisfying
\beql{ME18}
\la^2 - \la + 2 =0 ~.
\eeq
The first 32 terms in the 2-adic expansion \eqn{ME17} of $\la$ are
\beql{ME50}
0110010111111001110011011000110 \dd ~.
\eeq
There is no pattern to these digits.
}
\end{exam}

Then the 2-adic code of length 7 and type $1^4$
with generator polynomial
$$X^3 + \la X^2 + ( \la -1)X-1$$
is the 2-adic lift of the familiar binary $[7,4]$ Hamming code.
The generator polynomials for the versions of this code over
$\ZZ_2$, $\ZZ_4 , \dd$ are:
\beql{ME20}
\begin{array}{rll}
\ZZ_2 & : & X^3 + X+1 \\ [+.05in]
\ZZ_4 & : & X^3 + 2X^2 + X-1 \\ [+.05in]
\ZZ_8 & : & X^3 -2X^2 -3X -1 \\ [+.05in]
\ZZ_{16} & : & X^3 + 6X^2 + 5X -1 \\ [+.05in]
\ZZ_{32} & : & X^3 + 6X^2 + 5X-1 \\ [+.05in]
& ~ & \cdots
\end{array}
\eeq
(The coefficients can be read off \eqn{ME17}.)
By appending a 1 to the generating vectors of these codes, we obtain a sequence $\sH_2$, $\sH_4$, $\sH_8 , \dd , \sH_\infty$ of self-dual
codes.
In particular,
\beql{ME21}
\begin{array}{|ccccccc|c|}
\multicolumn{1}{c}{0} & 1 & 2 & 3 & 4 & 5 & \multicolumn{1}{c}{6} & \multicolumn{1}{c}{\infty} \\ \hline
1 & \la & \la -1 & -1 & 0 & 0 & 0 & 1 \\
0 & 1 & \la & \la-1 & -1 & 0 & 0 & 1 \\
0 & 0 & 1 & \la & \la -1 & -1 & 0 & 1 \\
0 & 0 & 0 & 1 & \la & \la -1 & -1 & 1 \\ \hline
\end{array}
\eeq
is the generator matrix for a self-dual 2-adic code $\sH_\infty$ of length 8 and type $1^4$ that we call the
{\em 2-adic Hamming code}.
This is in some sense the smallest interesting 2-adic code.

The $\ZZ_2$ version of this code, $\sH_2$, is the $[8,4]$ Hamming code, and the $\ZZ_4$ version, $\sH_4$,
is the {\em octacode}, studied in
\cite{CS93}, \cite{CS293}, \cite{Fo94}, \cite{HKCSS},
and equivalent to the binary nonlinear Nordstrom-Robinson code.

The minimal Hamming and Lee distances of these codes are as follows:
$$
\begin{array}{cccccccc}
~ & \sH_2 & \sH_4 & \sH_8 & \sH_{16} & \sH_{32} & \sH_{64} & \cdots \\
\mbox{Hamming} & 4 & 4 & 4 & 4 & 4 & 4 & \cdots \\
\mbox{Lee} & 4 & 6 & 8 & 12 & 14 & 18 & \cdots
\end{array}
$$
The minimal Hamming distance of $\sH_{2^a}$ for $1 \le a < \infty$ is always 4,
since the codeword obtained by multiplying any of the generators by
$2^{a-1}$ has Hamming weight 4.
However it follows from Theorem~\ref{Th8} below that the
2-adic
Hamming code $\sH_\infty$ has minimal Hamming distance 5, and is an MDS code.

On the other hand the sequence of Lee distances of these codes,
$4,6,8,12,14,18, \dd$, approaches infinity as $a \to \infty$.
Unfortunately it appears that this sequence does not converge
2-adically, so one obvious definition of the minimal Lee distance
of $\sH_\infty$ fails.
Even the Lee weight of the projections of the integer $\la$ onto
$\ZZ_{2^m}$ do not converge 2-adically as $m\to \infty$.
For let $\la = \sum\limits_{i=0}^\infty \la_i 2^i$ (the $\la_i$ are given in
\eqn{ME17}, \eqn{ME20}), so the projection onto
$\ZZ_{2^m}$ is $\alpha_m = \sum\limits_{i=0}^{m-1} \la_i 2^i$,
$m \ge 1$.
The Lee weight of $\alpha_m$ is $w_m = \min \{ \alpha_m , 2^m - \alpha_m \}$, and one can show that
$$w_m = (1-2 \la_{m-1}) \alpha_{m-1} + \la_{m-1} 2^{m-1}~,~~~
m \ge 2 ~.$$
This shows that $\{w_1 , w_2 , \dd \} = \{ 0,2,2,6,6,26, \dd \}$ does not
converge 2-adically.

There are several other natural ways to define the minimal distance
of this code, but none are completely satisfactory.
This is a question that requires further investigation.

The automorphism group of $\sH_\infty$ contains
operations corresponding to $x \mapsto x+1$, $x \mapsto 2x$ and $x \mapsto -1/x$, namely the monomials
$$
\begin{array}{l}
(0,1,2,3,4,5,6) ( \infty ) ~, \\
(0) (1,2,4) (3,6,5) ( \infty)~, \\
(0, \infty) (1,6) (2,3) (4,5) \mbox{~\& negate~} 0,1,2,4~,
\end{array}
$$
which generate the central product $\ZZ_2 . PSL_2 (7)$, as well as all scalar matrices $uI$, $u= \mbox{unit}$ in $\ZZ_{2^\infty}$.
Then the full projective automorphism group of $\sH_\infty$ is
$PSL_2 (7)$, of order 168.
\begin{exam}
\label{ex2}
{\bf The 2-adic Golay code of length 24}.
{\rm The binary Golay code can be lifted in a similar way.
The factorization of $X^{23} -1$ over $\ZZ_{2^\infty}$ is
$$X^{23} -1 = (X-1) \pi_\infty^{(1)} (X) \pi_\infty^{(2)} (X)~,$$
where
\begin{eqnarray}
\pi_\infty^{(1)} (X) & = &
X^{11} + \nu X^{10} + (\nu -3)X^9 - 4X^8 - (\nu +3) X^7 \nonumber \\
& ~ &~~- (2 \nu +1) X^6 - (2 \nu -3 ) X^5 - ( \nu -4 ) X^4 + 4X^3 \nonumber \\
& ~ & ~~+ (\nu+2)X^2 + ( \nu -1) X-1 ~,
\end{eqnarray}
\beql{ME23}
\nu = 0+2+8+32+64+128+ \cdots
\eeq
is a 2-adic number satisfying
\beql{ME24}
\nu^2 - \nu + 6 =0 ~,
\eeq
and $\pi_\infty^{(2)} (X)$ is the reciprocal polynomial to
$\pi_\infty^{(1)} (X)$.
The first 32 terms in the 2-adic expansion \eqn{ME23} are
$$0101011110010010110010000110000 \dd ~.$$
Then the cyclic code generated by $\pi_\infty^{(1)} (X)$,
extended by appending a 1 to the generators, is a self-dual
2-adic code $\sG_\infty$ of length 24 and type $1^{12}$, the
{\em 2-adic Golay code}.
The full projective automorphism group of $\sG_{\infty}$ is $PSL_2 (23)$.
}
\end{exam}

The projection on $\ZZ_2$ of $\sG_\infty$ is the binary Golay
code $\sG_2$ of length 24 and minimal Hamming distance 8,
and in fact every projection $\sG_{2^a}$ of this code onto
$\ZZ_{2^a}$ for finite $a$ has minimal Hamming distance 8.
However
it follows from Theorem~\ref{Th8} that the
2-adic Golay code $\sG_\infty$ has minimal Hamming distance 13,
and is an MDS code.

As in the previous example, the $\ZZ_4$ version of this code,
$\sG_4$, is especially interesting.
Bonnecaze and Sol\'{e} \cite{BS93} have shown that by applying
Construction~A to this code, i.e. by taking all vectors in $\ZZ^{24}$ which project onto $\sG_4$ modulo 4, one obtains the Leech lattice.
This is one of the simplest constructions known for this lattice (cf. \cite{CS93}).
\begin{exam}
\label{ex3}
{\bf The 3-adic Golay code of length 12}.
{\rm We lift the ternary Golay code in the same way,
using the irreducible divisor
$$X^5 + \theta X^4 - X^3 + X^2 + ( \theta -1) X-1$$
of $X^{11} -1$ over $\ZZ_{3^\infty}$, where
$$\theta = 0 + 3 + 9 + 2.27 + 2.81 + \cdots$$
is a 3-adic number satisfying
\beql{ME25}
\theta^2 - \theta +3 =0 ~.
\eeq
By appending a 1 to each generator we obtain a self-dual
3-adic code $\sT_\infty$ of length 12
and type $1^6$, the {\em 3-adic Golay code}.
This has minimal Hamming distance 7 and is an MDS
code.
Its full projective automorphism group is $PSL_2 (11)$.
}
\end{exam}
\begin{exam}
\label{ex4}
{\bf Binary quadratic residue codes}.
{\rm
Examples~\ref{ex1} and \ref{ex2} may be generalized as follows.
Let $n$ be a prime of the form $8m-1$, so that $X^n -1$ factorizes
over $\ZZ_2$ into $(X-1) \pi_2^{(1)} (X) \pi_2^{(2)} (X)$,
where all the factors are irreducible, with a corresponding
factorization $(X-1) \pi_\infty^{(1)} (X) \pi_\infty^{(2)} (X)$ over
$\ZZ_{2^\infty}$.
Let $Q$ and $N$ denote the nonzero quadratic residues and nonresidues
modulo $n$, and set
$$f_Q (X) = \sum_{i \in Q} X^i~,~~~
f_N (X) = \sum_{i \in N} X^i ~.$$
Then as in the binary case there are two inequivalent 2-adic
quadratic residue codes of length $n$.
}
\end{exam}
\begin{theo}
\label{Th9}
The two quadratic residue codes of prime length $n=8m-1$ over $\ZZ_{2^\infty}$ have
generator polynomials $\pi_\infty^{(1)}$ and $(X-1) \pi_\infty^{(1)} (X)$,
and idempotents
$$\alpha 1 + \beta f_Q (X) + \gamma f_N (X) ~,$$
where the coefficients $\alpha , \beta , \gamma$ are
the 2-adic numbers
$$\af = \frac{n+1}{2n} ,~~\beta = \frac{1+ \sqrt{-n}}{2n} , ~~
\gamma = \frac{1- \sqrt{-n}}{2n}
$$
for the first code, and
$$
\af = \frac{n-1}{2n} ,~~
\beta = \frac{-1 + \sqrt{-n}}{2n} , ~~
\gamma = \frac{-1- \sqrt{-n}}{2n}
$$
for the second code.
By appending $\sqrt{\frac{-1}{n}}$ to each generator of the first
code we obtain a self-dual code of length $n+1$ and type
$1^{(n+1)/2}$.
\end{theo}

We omit the straightforward proof, which includes
the verification that when $n=7$ and 23 the codes
generated by $\pi_\infty^{(1)} (X)$ coincide with those constructed
in Examples~\ref{ex1} and \ref{ex2}.
The full projective automorphism group of the self-dual
code of length $n+1$ is $PSL_2 (n)$.
\begin{theo}
\label{Th8}
The self-dual extended quadratic residue code of length $n+1$ described in Theorem~\ref{Th9} has minimal Hamming distance $(n+3)/2$, and is an
MDS code.
\end{theo}
\pf
It follows from Blahut \cite{Bl91} that this code consists of all vectors
$(c_0 , c_1 , \dd , c_{n-1} , c_\infty ) \in \ZZ_{2^\infty}^{n+1}$ that satisfy
$$
\sqrt{\frac{-1}{n}} \, \sum_{j=0}^{n-1} c_j + c_\infty =0 ~,
$$
$$
\sum_{j=0}^{n-1} c_j  \, \xi^{jq} =0 ,~~~q \in Q ~,
$$
where $\xi = e^{2 \pi i/n}$.
The usual Vandermonde argument then shows that this is an MDS code. \hfill $\blacksquare$
\begin{exam}
\label{ex5}
{\bf Cyclic codes of length 7 over $\ZZ_4$ and $\ZZ_{2^\infty}$}.
{\rm As an illustration of the structure theorems of Section~3
(and also because one of them is the octacode) we enumerate the cyclic
codes of length 7 over $\ZZ_4$.
We factorize $X^7 -1$ over $\ZZ_4$ from \eqn{ME16}, obtaining
\beql{ME40}
(X-1) (X^3 +2X+X-1) (X^3-X^2 +2X-1) = f_0 f_1 f_2
\eeq
(say).
The nontrivial prime ideals are,
from Theorem~\ref{Th2},
$$P_0 = (f_0 , 2), ~~P_1 = (f_1 ,2) ,~~P_2 = (f_2 ,2) ~,$$
and the other primary ideals are
$$P_0^2 = (f_0) , ~~P_1^2 = (f_1) , ~~ P_2^2 = (f_2 )~.$$
There are 27 codes, by Theorem~\ref{Th5}, and they are displayed in Table~1
(except that we have omitted codes $4,6, \dd , 27$, which
are equivalent to codes $3,5, \dd, 26$ under the
symmetry interchanging $f_1$ and $f_2$).
The fourth column gives the canonical forms for these codes
as described in Theorems~\ref{Th5} and \ref{Th6}.
}
\end{exam}

In Examples~\ref{ex1}--\ref{ex4} we extended the codes to length
$n+1$ by appending a symbol that made them self-dual.
For the codes in Table~1 it is more appropriate to append a zero-sum check symbol.
The two extensions agree in the case of the octacode, which is number 12.
The second column gives representative generators for the cyclic code
(with the extending symbol in parentheses).
The last column gives the minimal Lee distance $d$ of
the cyclic code (and the minimal distance $d^\ast$ of the extended
code in parentheses).
\begin{table}[htb]
\caption{Cyclic (and extended cyclic) codes of length 7 over $\ZZ_4$.
Number 12 is the octacode.}
$$
\begin{array}{rlccc}
\# & \mbox{generators} & \mbox{type} & \mbox{ideal} & d(d^\ast) \\ \hline
1 & 0000000 (0) & 1^0 & 0 = P_0^2 P_1^2 P_2^2 & \mbox{-(-)} \\
2 & 2222222 (2) & 2^1 & (2 f_1 f_2 ) = P_0 P_1^2 P_2^2 & 14 (16) \\
3 & 2220200 (0) & 2^3 & (2 f_0 f_1) = P_0^2 P_1^2 P_2 & 8(8) \\
5 & 2022000 (2) & 2^4 & (2 f_1) = P_0 P_1^2 P_2 & 6(8) \\
7 & 2200000 (0) & 2^6 & (2f_0) = P_0^2 P_1 P_2 & 4(4) \\
8 & 2000000 (2) & 2^7 & (2) = P_0 P_1 P_2 & 2(4) \\
9 & 1000000 (1) & 1^7 & (1)=1 & 1(2) \\
10 & 1300000 (0) & 1^6 & (f_0) = P_0^2 & 2(2) \\
11 & 1300000 (0), 2000000 (0) & 1^6 2^1 & (f_0,2) =P_0 & 2(2) \\
12 & 1213000 (1) & 1^4 & (f_1) = P_1^2 & 4(6) \\
14 & 1213000 (1), 2000000 (0) & 1^4 2^3 & (f_1,2) = P_1 & 2(4) \\
16 & 1132100 (0) & 1^3 & (f_0f_1) = P_0^2P_1^2 & 6(6) \\
18 & 1132100 (0), 2000000 (0) & 1^3 2^4 & (f_0 f_1,2) = P_0 P_1 & 2(4) \\
20 & 1132100 (0), 2200000 (0) & 1^3 2^3 & (f_0 f_1, 2 f_0) = P_0^2 P_1 & 4(4) \\
22 & 1132100 (0), 2022000 (2) & 1^3 2^1 & (f_0 f_1, 2 f_1) = P_0 P_1^2 & 4(6) \\
24 & 1111111 (1) & 1^1 & (f_1 f_2) = P_1^2 P_2^2 & 7(8) \\
25 & 1111111 (1), 2000000 (0) & 1^1 2^6 & (f_1 f_2,2) = P_1 P_2 & 2(4) \\
26 & 1111111 (1), 2022000 (0) & 1^1 2^3 & (f_1 f_2,2f_1) = P_1^2 P_2 & 6(8)
\end{array}
$$
\end{table}

It is easy to extend this table to obtain a list of all possible
types of cyclic codes over length $n$ over
$\ZZ_q$, $q=p^a$, $1 \le a \le \infty$, for any
prime $p$ such that $X^n -1$ factorizes modulo $p$ into three irreducible factors, as in \eqn{ME40}.
It follows from Theorem~\ref{Th6} that there are 24 types of such codes, namely
$$(p^{m_0} g_0) ,~~(p^{m_0} g_0 , p^{m_1}) ~,$$
where
$g_0 \in \{ f_0, f_1, f_2, f_0 f_1, f_0 f_2, f_1 f_2 \}$,
and
$$(p^{m_0} g_0 , p^{m_1} g_1 ) , ~~
(p^{m_0} g_0 , p^{m_1} g_1 , p^{m_2} ) ~,$$
where $g_0 \in \{ f_0 f_1 , f_0 f_2 , f_1 f_2 \}$,
$g_1 | g_0$, and
$$0 \le m_0 < m_1 < m_2 ~.$$
Similar enumerations can be obtained for any $n$,
once the factorization of $X^n -1$ is known.
\subsection*{Acknowledgements}
\hsp
We thank Mira Bernstein, Joe Buhler and especially John Conway
for helpful conversations,
and Christine Chang for assistance in tabulating
cyclic codes over $\ZZ_4$.

\end{document}